\documentclass[12pt,leqno]{article} 
\usepackage{graphics}

\newcommand{\beqa}{\begin{eqnarray}}
\newcommand{\eeqa}{\end{eqnarray}}

\newcommand{\beq}{\begin{equation}}
\newcommand{\eeq}{\end{equation}}
\newcommand{\lbl}{\label}
\newcommand{\s}{\; \;}

\title{A remark on Pinney's equation}

\author{
Philip Korman   \\ 
Department of Mathematical Sciences \\ 
University of Cincinnati \\ 
Cincinnati Ohio 45221-0025 \\
}

\date{}

\begin{document}

\maketitle
\begin{abstract} 
We show that Pinney's equation \cite{P} with a constant coefficient can be reduced to its linear part by a simple change of variables. Also, Pinney's original solution is simplified slightly.  
 \end{abstract}

\begin{flushleft}
\end{flushleft}

\begin{flushleft}
\end{flushleft}

In 1950 Edmund Pinney published a very influential paper \cite{P}, which was less than half a page long. That paper provided a general solution of the nonlinear differential equation
\beq
\lbl{1}
y''+a(x)y+\frac{c}{y^3}=0 \,, \s y(x_0)=q \ne 0\,, \; y'(x_0)=p \,,
\eeq
with a given function $a(x)$ and a constant $c \ne 0$. Namely, the solution is
\beq
\lbl{2}
y(x)=\pm \sqrt{u^2(x)-cv^2(x)} \,,
\eeq
where $u(x)$ and $v(x)$ are the solutions of the linear equation 
\beq
\lbl{3}
y''+a(x)y=0 \,,
\eeq
for which $u(x_0)=q$, $u'(x_0)=p$, and $v(x_0)=0$, $v'(x_0)=\frac{1}{q}$. One takes ``plus" in (\ref{2}) if $q>0$, and ``minus" if $q<0$. Clearly, $u(x)$ and $v(x)$ form a fundamental set, and by Liouville's formula their Wronskian
\[
u'(x)v(x)-u(x)v'(x)=1 \,, \s \mbox{for all $x$}\,.
\]
A substitution of $y=\sqrt{u^2(x)-cv^2(x)}$ into (\ref{1}) gives
\[
y''+a(x)y+\frac{c}{y^3}=-c \, \frac{\left[u'(x)v(x)-u(x)v'(x)\right]^2-1}{\left[u^2(x)-cv^2(x)\right]^{\frac{3}{2}}} =0 \,.
\]
If $c<0$, the solution is valid for all $x$, while for $c>0$ some singular points are possible.
\medskip

The nonlinear  equation equation (\ref{1}) possessing a {\em general solution} is very special, and it attracted a lot of attention (there are currently 92 MathSciNet and 543 Google Scholar citations). It turns out that this equation was considered back in 1880 by Ermakov \cite{E}.
\medskip

Our remark is that in case of constant $a(x)=a_0$, Pinney's equation becomes linear for $z(x)=y^2(x)$. Indeed, we multiply the equation 
\beq
\lbl{4}
y''+a_0y+\frac{c}{y^3}=0
\eeq
by $y'$, and integrate to get
\beq
\lbl{5}
{y'}^2+a_0y^2-c{y}^{-2}=p^2+a_0q^2-c\frac{1}{q^2} \,.
\eeq
Now multiply the same equation by $y$:
\beq
\lbl{6}
yy''+a_0y^2+cy^{-2}=0 \,,
\eeq
and set $z=y^2$. Since $yy''=\frac{1}{2}z''-{y'}^2$, by using (\ref{5}), one transforms (\ref{6}) to 
\[
z''+4a_0z=2 \left(p^2+a_0q^2-c\frac{1}{q^2}\right) \,, \s z(0)=q^2 \,, \; z'(0)=2pq \,.
\]

\end{document}